\documentclass[11pt]{article}
\usepackage{amsmath,amssymb,amsthm,xcolor}
\usepackage[T1]{fontenc}

\newcommand{\eq}[1]{(\ref{eq:#1})}
\newcommand{\thm}[1]{Theorem~\ref{thm:#1}}

\newcommand{\rem}[1]{Remark~\ref{rem:#1}}
\newcommand{\cor}[1]{Corollary~\ref{cor:#1}}
\newcommand{\defi}[1]{Definition~\ref{def:#1}}
\newcommand{\prop}[1]{Proposition~\ref{prop:#1}}
\newcommand{\geo}{\text{\rm geo}}
\newcommand{\ggeo}{\text{\rm g-geo}}
\newcommand{\gexp}{\text{\rm g-exp}}
\newtheorem{Definition}{Definition}
\newtheorem{Theorem}{Theorem}

\newtheorem{Proposition}{Proposition}
\newtheorem{Corollary}{Corollary}
\newtheorem{Remark}{Remark}
\newtheorem{Example}{Example}

\begin{document}
\title{On independence of time and cause}

\author{Offer Kella,
\thanks{Department of Statistics; The Hebrew University of Jerusalem; Jerusalem 9190501; Israel. (\tt{offer.kella@gmail.com})} \thanks{Supported in part by the Vigevani Chair in Statistics.}
}

\date{\today}

\maketitle

\begin{abstract}
For two independent, almost surely finite random variables, independence of their minimum (time) and the events that either one of them is equal to the minimum (cause) is completely characterized. It is shown that, other than for trivial cases where, almost surely, one random variable is greater than or equal to the other, this happens if and only if both random variables are distributed like the same strictly increasing function of  two independent random variables, where either both are exponentially distributed or both are geometrically distributed. This is then generalized to the multivariate case.
\end{abstract}


\bigskip
\noindent {\bf AMS Subject Classification (MSC2020):} Primary 62E10; Secondary 60E05, 90B25.

\section{Introduction}

Let $X,Y$ be independent almost surely finite random variables. As is common in {\em e.g.,} queueing, reliability and discrete and continuous time Markov chains, we may think of $X,Y$ as competing time epochs of a certain event ({\em e.g.,} arrival, failure or transition) relative the current time and in general can take on any real value, not necessarily nonnegative. Let us refer to $\min(X,Y)$ as (event) {\em time} and to $\{X=\min(X,Y)\}$ and $\{Y=\min(X,Y)\}$ as (event) {\em causes}. When (the $\sigma$-field generated by) {\em time} is independent of (that generated by) the two {\em causes} are independent, we will say that {\em time and cause are independent}.

It is a simple undergraduate level exercise to show that if $X,Y$  are either both exponentially distributed or both geometrically distributed, then time and cause are independent. In the exponential case $P(Y=X)=0$ while in the geometric case $P(Y=X)>0$. Also, time and cause are independent for $X,Y$ if and only if they are independent for $g(X),g(Y)$ for any strictly increasing function $g$. Let us refer to the distribution of $g(X)$ where $X$ has an exponential (resp., geometric) distribution as
a $g$-exponential (resp., $g$-geometric) distribution.

There are trivial instances when time and cause must be independent that do not belong to the setup of the previous paragraph. Namely, when for some $c$ either $P(X<c)=P(Y\ge c)=1$ or $P(X\le c)=P(Y>c)=1$ or $P(X=c)=P(Y\ge c)=1$ (similarly, when $X,Y$ are interchanged).

The purpose of this paper is argue that the previous two paragraphs exhaust all the cases where time and cause can be independent. The results regarding the $g$-exponential and $g$-geometric cases are then generalized to the multivariate case.

In \cite{s65} it was shown that if $X,Y$ are independent random variables with continuous distribution functions $F_X,F_Y$, respectively and $P(X>Y)$ and $P(X<Y)$ are both positive, then $\min(X,Y)$ and $\{Y>X\}$ are independent if and only if $F_Y(t)=1-(1-F_X(t))^\lambda$ for some $\lambda>0$ and all $t$. This is equivalent to the statement that the corresponding cumulative hazard functions are proportional.
The case where $F_X,F_Y$ are absolutely continuous was discussed in \cite{a58,a63} and again in \cite{n70}, the latter citing only \cite{a63}, apparently unaware of \cite{a58,s65}. Some discussion where $X,Y$ are not independent and have a jointly absolutely continuous distribution can be found in \cite{kp91}.

In relation to proportional hazards it should be mentioned that proportional hazards models have an old and substantial history in statistics. So common that it even has an entry in Wikipedia (Proportional hazards model) with a substantial (but certainly far from exhaustive) reference list. Related statistical models (in particular the Cox proportional hazards model) appears in textbooks and taught in various statistics courses. Since this is not the focus of this paper, I hope that the reader forgives me for not mentioning specific references which can be quickly and readily found via a simple internet search.

My own initial motivation for this study was a setup where $X,Y$ are independent, $X$ has an exponential distribution and $Y$ satisfies $P(Y>0)>0$. It was important to know whether it is true that time and cause can be independent only when $Y$ is exponentially distributed. If it is assumed that the distribution of $Y$ is continuous, then this would follow from \cite{s65}. However, this turns out to be true without this a-priori assumption. To be honest, for this particular purpose, \cite{s65} is sufficient once one separately argues that under these assumptions $Y$ must also have a continuous distribution, which is not hard, but since $P(X=Y)=0$, this is an immediate consequence of \thm{main} below (see \rem{example}).

In what follows, Section~\ref{setup} contains the main results, while Section~\ref{discussion} contains some remarks.

\section{Setup and main results}\label{setup}
In what follows, {\em cumulative distribution function}, {\em almost surely} and {\em if and only if} are abbreviated to {\em cdf}, {\em a.s.} and {\em iff} respectively. Denote,
$x\wedge y=\min(x,y)$. For a random variable $X$, $F_X(x)=P(X\le x)$ denotes the associated cdf. All random variables considered here are a.s. finite. All cdf's are of a.s. finite random variables. When $F$ is a cdf, we denote $\bar F(\cdot)=1-F(\cdot)$. The notation `$\sim$' abbreviates {\em distributed} and {\em distributed like}. Finally, a random variable $X$ is independent of an event $A$ means that $P(\{X\le t\}\cap A)=F_X(t)P(A)$ for all $t$ and, equivalently, that $P(\{X\in B\}\cap A)=P(X\in B)P(A)$ for all Borel $B$.

From here on, $g$ will be a (resp., strictly increasing) function with domain $(0,\infty)$. If $h$ is a (resp., strictly increasing) function with domain $\{1,2,\ldots\}$, then there always is a (resp., strictly increasing) function $g$ with domain $(0,\infty)$ which coincides with $h$ on $\{1,2,\ldots\}$. By this we avoid the nuisance of specifying the domain of $g$ every time and in each case (exponential or geometric).

\begin{Definition}
When $X=X\wedge Y$ (resp., $Y=X\wedge Y$) we will say that $X$ (resp., $Y$) is a cause.
\end{Definition}

Note that when $X=Y$ then both $X$ and $Y$ are causes. It is easily checked that an event or a random variable is independent of the two events $\{X=X\wedge Y\}$ and $\{Y=X\wedge Y\}$ iff it is independent of the three events $\{X<Y\}$, $\{X=Y\}$ and $\{X>Y\}$. Therefore, we find it more convenient to define independence of time and cause as follows.
\begin{Definition}\label{def:timecause}
For random variables $X,Y$, time and cause are independent means that the random variable $X\wedge Y$ is independent of the events $\{X>Y\}$, $\{X<Y\}$ and $\{X=Y\}$.
\end{Definition}
Clearly, if $X\wedge Y$ is independent of any two of the three events in Definition~\ref{def:timecause}, then it is also independent of the third. In particular, if one of the events has probability zero ({\em e.g.,} when one of the distributions is known to be continuous, and then $P(X=Y)=0$), then it is enough to assume that $X\wedge Y$ is independent of one of the other two events.

\begin{Definition}
$X$ has a $g$-exponential (resp., $g$-geometric) distribution if $X\sim g(V)$ and $V$ has an exponential (resp., geometric) distribution. We write $X\sim \gexp(\lambda)$ (resp., $\ggeo(\alpha)$) when we want to emphasize that $V\sim\exp(\lambda)$ (resp., $\geo(\alpha)$) for some specific $\lambda>0$ (resp., $\alpha\in (0,1)$.
\end{Definition}
When writing that {\em $X_1,\ldots,X_n$ have $g$-exponential distributions,} we mean that $X_i\sim \gexp(\lambda_i)$ for $i=1,\ldots,n$ and some (possibly different) $\lambda_1,\ldots,\lambda_n$, while $g$ is the same for all $i$. Similarly for $g$-geometric.

Let us first get the trivial cases out of the way.
\begin{Proposition}\label{prop:main}
Assume $P(X>Y)=0$. If $P(X=Y)=0$, then time and cause are independent.  If $P(X=Y)>0$ then time and cause are independent
iff, for some $c$, $P(X=c)=P(Y\ge c)=1$.
Similarly, when $X,Y$ are interchanged.
\end{Proposition}
Note that when $P(X=Y)=1$, then $P(X=c)=P(Y=c)=1$ for some $c$, a trivial case which is already included in \prop{main}. Also note that $P(X>Y)=P(X=Y)=0$ (so $P(X<Y)=1$), iff for some $c$, either $P(X\le c)=P(Y>c)=1$ or $P(X<c)=P(Y\ge c)=1$.

\begin{proof}[\bf Proof:]
If $P(X=Y)=0$ then $P(X<Y)=1$. In this case, all the three events from \defi{timecause} have probability zero or one and are thus independent of any random variable. In particular $X\wedge Y$.

If $P(X=Y)>0$, then for some $c$, $P(X\le c)=P(Y\ge c)=1$ and $P(X=c),P(Y=c)>0$. Denoting $Z=X\wedge Y$, we clearly have that $P(X=Y)=P(X=c)P(Y=c)$, while 
\begin{equation}
P(Z=c)=P(X=c)P(Y\ge c)=P(X=c)\,.\end{equation} 
When $Z$ is independent of $\{X=Y\}$ we must have that
\begin{align}
P(X=c)P(Y=c)&=P(Z=c,X=Y)=P(Z=c)P(X=Y)\nonumber\\
&=P(X=c)^2P(Y=c)\,,
\end{align}
which implies that $P(X=c)=1$. Therefore $P(Z=c)=1$ and thus $Z$ is independent of any event. In particular the events from \defi{timecause}.
\end{proof}

Before stating our main result let us recall the notion of a {\em generalized inverse} ({\em e.g.,} \cite{eh13}) of a cdf $F$, defined as follows for  every $u\in(0,1)$:
\begin{equation}\label{eq:finverse}
F^{-1}(u)=\inf\{x|F(x)\ge u\}\,.
\end{equation}
It is possible to define $F^{-1}$ for $u=0$ and $u=1$ but this will not be needed here. The following are some necessary (well known) properties.
\begin{description}
\item{(i)} If $U\sim\text{Uniform}(0,1)$ then $F^{-1}(U)\sim F$.
\item{(ii)} If $F$ is continuous and $X\sim F$, then $F(X)\sim \text{Uniform}(0,1)$.
\item{(iii)} If $F$ is continuous then $F^{-1}$ is strictly increasing on $(0,1)$.
\item{(iv)} If $F$ is continuous then $F(F^{-1}(u))=u$ for all $u\in(0,1)$.
\end{description}
From (ii) it follows that when $F$ is continuous then, with $H_F(x)=-\log\bar F(x)$ when $F(x)<1$ (cumulative hazard function), $H_F(X)\sim\exp(1)$. If we denote
\begin{equation}\label{eq:hinverse}
H^{-1}_F(y)=\inf\{x|\,H_F(x)\ge y\}\,,
\end{equation}
for $y>0$, then $H^{-1}_F(y)=F^{-1}(1-e^{-y})$ and thus, if $V\sim\exp(1)$ then $1-e^{-V}\sim U\sim \text{Uniform}(0,1)$, so so that $H^{-1}_F(V)\sim F^{-1}(U)\sim F$. When $F$ is continuous then, from (iii), $H^{-1}_F$ is strictly increasing on $(0,\infty)$. For $\lambda>0$ and $F$ continuous, denote $F_\lambda(\cdot)=1-\bar F^\lambda(\cdot)$. Then it is easy to check that $F_\lambda^{-1}(u)=F^{-1}(1-(1-u)^{1/\lambda})$, so that $H^{-1}_{F_\lambda}(y)=F^{-1}(1-e^{-y/\lambda})=H_F^{-1}(y/\lambda)$. If $V\sim\exp(1)$ then $V/\lambda\sim\exp(\lambda)$ so that $H_F^{-1}(V/\lambda)=H_{F_\lambda}^{-1}(V)\sim F_\lambda$. Thus, if $F$ is continuous, $X\sim F_\lambda$ and $Y\sim F_\mu$, then, with $g=H^{-1}_F$, $X\sim g(V)$ and $Y\sim g(W)$ where $V\sim\exp(\lambda)$, $W\sim\exp(\mu)$ and $g$ is strictly increasing. This implies that for a continuous cdf
\begin{equation}\label{eq:CF}
\mathcal{C}_F=\{F_\lambda(\cdot)|\,\lambda>0\}=\{F_{g(V)}|\,V\sim \exp(\lambda),\,\lambda>0\}\,,
\end{equation}
where either $F$ is given and on the right $g=H^{-1}_F$ is strictly increasing, or $g$ is given, strictly increasing, and on the left $F=F_{g(V)}$ for $V\sim\exp(1)$. Note that when $F_V$ is continuous and $g$ is strictly increasing, then $F_{g(V)}$ is also continuous. Also note that $\mathcal{C}_F=\mathcal{C}_{F_\mu}$ for all $\mu>0$, which is equivalent to replacing $V\sim\exp(\lambda)$ by $V\sim\exp(\lambda\mu)$ on the right hand side of~\eq{CF}.

\begin{Remark}\label{rem:expgeo}\rm
Assume that $X,Y$ are distributed like $g(V),g(W)$ where $V,W$ are independent and $g$ is strictly increasing. Then $Z\equiv X\wedge Y\sim g(V\wedge W)$, $P(X>Y)=P(V>W)$, $P(X=Y)=P(V=W)$ and $P(X<Y)=P(V<W)$.

Therefore, when $X\sim \gexp(\lambda)$ and $Y\sim \gexp(\mu)$, then, with $\nu=\lambda+\mu$,
\begin{align}
P(X=X\wedge Y)&=P(X<Y)=\frac{\lambda}{\nu}\nonumber\\
P(Y=X\wedge Y)&=P(Y>X)=\frac{\mu}{\nu}\\
Z&\sim \gexp(\nu)\,.\nonumber
\end{align} 
When $X\sim \ggeo(\alpha)$ and $Y\sim \ggeo(\beta)$, then, with $\gamma=1-(1-\alpha)(1-\beta)$,
\begin{align}
P(X<Y)&=\frac{\alpha(1-\beta)}{\gamma}\qquad P(X=Y)=\frac{\alpha\beta}{\gamma}\nonumber\\
P(X>Y)&=\frac{\beta(1-\alpha)}{\gamma}\qquad Z\sim \ggeo(\gamma)\,,
\end{align} 
so that
\begin{align}
P(X=X\wedge Y)&=P(X\le Y)=\frac{\alpha}{\gamma}\nonumber\\
P(Y=X\wedge Y)&=P(X\ge Y)=\frac{\beta}{\gamma}\,.
\end{align}
\end{Remark}

The following is the main result. In the proof we make some use of a specific notion of a cumulative hazard function which may be defined in more than one way. See~\cite{ks80} for a nice in depth analysis of relationships between hazard (and other) functions and probability distributions.
\begin{Theorem}\label{thm:main}
Let $X,Y$ be independent random variables with
\begin{equation}
p=P(X<Y)>0\,,\qquad q=P(X>Y)>0
\end{equation} and denote $r=P(X=Y)$. Then time and cause are independent iff, for some strictly increasing $g$, the distributions of $X,Y$, hence also that of $X\wedge Y$, are either $g$-exponential (when $r=0$) or $g$-geometric (when $r>0$).
\end{Theorem}

\begin{proof}[\bf Proof:]
The fact that the  $g$-exponential and $g$-geometric cases (for strictly increasing $g$) imply independence of time and cause is clear and was discussed in the introduction. Therefore, it remains to show that if time and cause are independent (and $p,q>0$) then the distributions of $X,Y$ are either both $g$-exponential or both $g$-geometric for some (common) strictly increasing $g$.

For a cdf $F$ we recall that $\bar F(\cdot)=1-F(\cdot)$ and denote $F(t-)=\lim_{s\uparrow t}F(s)$, $\Delta F(t)=F(t)-F(t-)$, $F^c(t)=F(t)-\sum_{s\le t}\Delta F(s)$ and finally $t_F=\sup\{t|\,F(t)<1\}$. It may be checked  ({\em e.g.,} \cite[Thm.~54, p.~41]{p04}), that for $t<t_F$, the cumulative hazard function $H_F(t)=-\log\bar F(t)$ satisfies
\begin{align}\label{eq:hazard}
H_F(t)&=\int_{-\infty}^t\frac{1}{\bar F(s)}dF^c(s)+\sum_{s\le t}(H_F(t)-H_F(t-))\nonumber\\
&=\int_{-\infty}^t\frac{1}{\bar F(s)}dF^c(s)+\sum_{s\le t}\log\left(1+\frac{\Delta F(t)}{\bar F(t)}\right)\,.
\end{align}
For a given random variable $X$, we will abreviate $t_{F_X}$ to $t_X$ and $H_{F_X}$ to $H_X$.

Denote (again) $Z=X\wedge Y$. Then $\bar F_Z(t)=\bar F_X(t)\bar F_Y(t)$. By the independence of time and cause we have that
\begin{align}\label{eq:delta1}
\Delta F_X(t)\bar F_Y(t)&=P(Z=t,X<Y)=p\Delta F_Z(t)\nonumber\\
\Delta F_Y(t)\bar F_X(t)&=P(Z=t,X>Y)=q\Delta F_Z(t)
\end{align}
and from 
\begin{equation}\label{eq:delta2}
\Delta F_Z(t)=\Delta F_X(t)\bar F_Y(t)+\bar F_X(t)\Delta F_Y(t)+\Delta F_X(t)\Delta F_Y(t)\,,
\end{equation} 
we also have that
\begin{equation}\label{eq:delta3}
\Delta F_X(t)\Delta F_Y(t)=r\Delta F_Z(t)\,. 
\end{equation}
From \eq{delta1} we immediately have that for every $t<t_Z=t_X\wedge t_Y$ (hence $\bar F_X(t),\,\bar F_Y(t)>0$), $\Delta F_X(t)$, $\Delta F_Y(t)$ and $\Delta F_Z(t)$ are either all positive or all zero. When they are positive we necessarily have from \eq{delta1}, \eq{delta2} and \eq{delta3} that, for $a=1+\frac{r}{p}$ and $b=1+\frac{r}{q}$,
\begin{align}\label{eq:ratio}
1+\frac{\Delta F_X(t)}{\bar F_X(t)}&=a\nonumber\\
1+\frac{\Delta F_X(t)}{\bar F_Y(t)}&=b\\
1+\frac{\Delta F_Z(t)}{\bar F_Z(t)}&=1+\frac{r}{q}+\frac{r}{p}+\frac{r^2}{pq}=ab\,.\nonumber
\end{align}
From \eq{ratio} and \eq{hazard} it now follows that, for any $t<t_z$, the number of discontinuities $N(t)$ on $(-\infty,t]$ is the same for $F_X$, $F_Y$, $F_Z$, is finite and
\begin{align}\label{eq:hazard1}
H_X(t)&=\int_{-\infty}^t\frac{1}{\bar F_X(s)}dF^c_X(s)+N(t)\log a\nonumber\\
H_Y(t)&=\int_{-\infty}^t\frac{1}{\bar F_Y(s)}dF^c_Y(s)+N(t)\log b\\
H_Z(t)&=\int_{-\infty}^t\frac{1}{\bar F_Z(s)}dF^c_Z(s)+N(t)\log ab\,.\nonumber
\end{align}
When $t_Z<\infty$ then $F_Z$ must be continuous at $t_Z$. To see this, assume that $\Delta Z(t_Z)>0$. Then from \eq{delta1} it follows that $\bar F_X(t_Z)>0$ and $\bar F_Y(t_Z)>0$. Thus, $\bar F_Z(t_Z)=\bar F_X(t_Z)\bar F_Y(t_Z)>0$, which contradicts the definition of $t_Z$.

Noting that $P(Z>t,\,X<Y)=P(t<X<Y)$ (and similarly when $X,Y$ are interchanged), then due to independence of time and cause, we have the following.
\begin{align}\label{eq:ind1}
P(Z>t)=\frac{P(t<X<Y)}{p}=\frac{P(t<Y<X)}{q}\,.
\end{align}
We observe from \eq{delta1} that for each $t<t_Z$,
\begin{align}
P(t<X<Y)&=\int_{(t,\infty)}\bar F_Y(t)dF_X(t)\nonumber\\
&=\int_{(t,\infty)}\bar F_Y(s)dF_X^c(s)+\sum_{s>t}\bar F_Y(s)\Delta F_X(s)\\
&=\int_{(t,\infty)}\bar F_Y(s)dF_X^c(s)+p\sum_{s>t}\Delta F_Z(s)\,,\nonumber
\end{align}
and similarly with $P(t<Y<X)$. Therefore, for $t<t_Z$ we have that
\begin{equation}
dF_Z^c(t)+\Delta F_Z(t)=\frac{1}{p}\bar F_Y(t)dF_X^c(t)+\Delta F_Z(t)=\frac{1}{q}\bar F_X(t)d F^c_Y(t)+\Delta F_Z(t)\,,
\end{equation}
so that
\begin{equation}
dF_Z^c(t)=\frac{1}{p}\bar F_Y(t)dF_X^c(t)=\frac{1}{q}\bar F_X(t)dF_Y^c(t)\,.
\end{equation}
Dividing by $\bar F_Z(t)=\bar F_X(t)\bar F_Y(t)$ gives
\begin{equation}
\frac{1}{\bar F_Z(t)}dF_Z^c(t)=\frac{1}{p}\frac{1}{\bar F_X(t)}dF_X^c(t)=\frac{1}{q}\frac{1}{\bar F_Y(t)}dF_Y^c(t)\,.
\end{equation}
Denoting 
\begin{equation}
F(t)=\begin{cases}
1-\exp\left(-\int_{(-\infty,t]}\frac{1}{\bar F_Z(s)}dF^c_Z(s)\right)&t<t_Z\\
1&t\ge t_Z\,,
\end{cases}
\end{equation} 
then $F$ is continuous, possibly with $F(t)=1$ for all $t$, when $F_Z^c(t)=0$ for all $t$ (the cdf of a random variable which is infinity a.s.). Recalling~\eq{hazard1}, we have, for $t<t_Z$ 
, that
\begin{align}\label{eq:FXYZ}
F_X(t)&=1-\bar F(t)^pb^{-N(t)}\nonumber\\
F_Y(t)&=1-\bar F(t)^qa^{-N(t)}\\
F_Z(t)&=1-\bar F(t)(ab)^{-N(t)}\,.\nonumber
\end{align}
Therefore, for $t<t_Z$,
\begin{equation}\label{eq:r}
\bar F(t)(ab)^{-N(t)}=\bar F_Z(t)=\bar F_X(t)\bar F_Y(t)=\bar F(t)^{1-r}(ab)^{-N(t)}\,.
\end{equation}
If $r=0$, then $a=b=1$ and we have that
\begin{equation}\label{eq:prop}
\bar F(t)=\bar F_Z(t)=\bar F_X(t)^{1/p}=\bar F_Y(t)^{1/q}\,,
\end{equation} 
for $t<t_Z$. Since $\bar F_Z(t)\to 0$ as $t\uparrow t_Z$ ($F_Z$ is continuous at $t_Z$), then we also have that $\bar F_X(t),\bar F_Y(t)\to 0$ and thus \eq{prop} is satisfied for every $t$. This implies that in \cite[Th.~1]{s65}, the condition that $F_X$ and $F_Y$ are continuous may be replaced with the weaker condition that $P(X=Y)=0$. From \eq{CF} it follows that in this case the distributions of $X,Y$ are $g$-exponential for $g=H_F^{-1}$, which, since $F$ is continuous, is strictly increasing.

On the other hand, if $r>0$, then since $\bar F$ is positive on $(-\infty,t_Z)$ and $N(t)<\infty$ for every $t<t_Z$, it follows from \eq{r} that necessarily $\bar F(t)=1$ for all $t<t_Z$. In this case, since $F_Z(\cdot)$ is continuous at $t_Z$ and $F_Z(t_Z)=0$, then we also must have, by \eq{FXYZ}, that $N(t)\to\infty$ as $t\uparrow t_Z$. For $i\ge 1$ let $g(i)=\inf\{t|\,N(t)=i\}$ be the discontinuities of $F_Z$. Then $g$ is clearly increasing on the positive integers and it follows from \eq{FXYZ} (with $\bar F(t)=1$ for all $t$) that $X\sim g(N_{1-a^{-1}})$ and $Y\sim g(N_{1-b^{-1}})$, where $N_\alpha\sim\geo(\alpha)$. The proof is complete.
\end{proof}

\prop{main} and \thm{main} give a complete characterization of when time and cause are independent for a pair of independent random variables $X,Y$ which, together with \rem{expgeo}, are summarized by the following Table~1, in which $p=P(X<Y)$, $q=P(X>Y)$ and $r=P(X=Y)$. In this table, $0$ and $+$ abbreviate $=0$ and $>0$, respectively. Also $c,\lambda,\mu,\alpha,\beta,g$ in Table~1 mean: {\bf for some} $c\in\mathbb{R}$, $\lambda,\mu>0$, $\alpha,\beta\in (0,1)$, strictly increasing $g$, respectively. Finally $\nu=\lambda+\mu$ and $\gamma=1-(1-\alpha)(1-\beta)$. Recall that $p+r=P(X=X\wedge Y)$ and $q+r=P(Y=X\wedge Y)$ are the probabilities that $X$ and $Y$, respectively, are causes. 

  \begin{center}
    \begin{tabular}{c|c|c|c|c|c|}\\ 
      $p$ & $q$ & $r$ & Conditions & $p+r$ & $q+r$\\
      \hline
      $+$ & $0$ & $0$ & None & 1 & 0\\
      $0$ & $+$ & $0$ & None & 0 & 1\\
      $0$ & $0$ & $+$ & None & 1 & 1\\
      $+$ & $0$ & $+$ & $P(X=c)=P(Y\ge c)=1$ & 1 &  $\Delta F_Y(c)$\\
      $0$ & $+$ & $+$ & $P(Y=c)=P(X\ge c)=1$ & $\Delta F_X(c)$ & 1\\
      $+$ & $+$ & $0$ & $X\sim \gexp(\lambda),\,Y\sim \gexp(\mu)$ & $\lambda/\nu$ &$\mu/\nu$ \\ 
      $+$ & $+$ & $+$ & $X\sim \ggeo(\alpha),\,Y\sim \ggeo(\beta)$ & $\alpha/\gamma$ & $\beta/\gamma$
      \end{tabular}
      \vskip .3 cm
      Table 1
  \end{center}
In lines $+|+|0$ and $+|+|+$ the distribution of $X\wedge Y$ is $\gexp(\nu)$ and $\ggeo(\gamma)$, respectively. Otherwise $X\wedge Y$ is either $X$ a.s. or $Y$ a.s. (or both).

Up until now we have dealt with independence of time and cause for two random variables. We now proceed to generalize these ideas to the multivariate case, which fortunately turns out to be quite straightforward.
\begin{Definition}\label{def:multi}
We will say that $X_1,\ldots,X_n$ satisfy that time and cause are independent when $\min_{1\le i\le n}X_i$ is independent of the events 
\begin{equation}
\{X_i>\min_{j\not=i}X_j\},\ \{X_i=\min_{j\not=i}X_j\},\ \{X_i<\min_{j\not=i}X_j\}\,.
\end{equation} 
for all $1\le i\le n$.
\end{Definition}
We note that in general for $n\ge 3$, the condition in \defi{multi} is stronger than the condition that $\min_{1\le i\le n}X_i$ is independent of $\{X_i=\min_{1\le j\le n}X_j\}$ for all $i$.
It is equivalent when $P(X_i=X_j)=0$ for all $i\not=j$.

\begin{Corollary}\label{cor:main}
Assume that $X_1,\ldots,X_n$ are independent and assume that $P(X_i>\min_{j\not=i}X_j)>0$ and $P(X_i<\min_{j\not=i}X_j)>0$ for all $1\le i\le n$.
Then time and cause are independent iff, for some strictly increasing $g$, $X_1,\ldots,X_n$ are all $g$-exponentially distributed or all $g$-geometrically distributed. 
\end{Corollary}

\begin{proof}[\bf Proof:]
As for the case $n=2$, in the $g$-exponential an $g$-geometric cases, time and cause are clearly independent. For the converse, denote $M=\min_{1\le i\le n}X_i$ and, for each $1\le i\le n$, $M_i=\min_{j\not=i}X_j$. Then $X_i\wedge M_i=M$. Thus, by \thm{main} and \rem{expgeo}, for any given $i$, $X_i,M_i,M$ are either  $g$-exponentially or $g$-geometrically distributed for some strictly increasing $g$. We recall that $g$ was determined from the distribution of $M$. Namely, for the $g$-exponential distribution we took $g=H_M^{-1}$ and for the $g$-geometric case $g(k)$ are the discontinuity points of $F_M$. Therefore, the same $g$ applies to all $i$ and either $X_i$ are $g$-exponentially distributed for all $i$ or they are $g$-geometrically distributed for all $i$.
\end{proof}

\section{Some discussion}\label{discussion}
In the following remarks we will now make some straightforward observations.
\begin{Remark}\label{rem:example}\rm
From \cor{main} and the notations in its proof, it also follows that, under its assumption, either $P(X_i=M_i)=0$ for all $i=\ldots n$ or $P(X_i=M_i)>0$ for all $i=1,\ldots,n$.
Also note that if $P(X_i=M_i)=0$ (resp., $>0$) for any $i$, then $M$ must have a $g$-exponential (resp., $g$-geometric) distribution for some increasing $g$. Therefore, so do $X_i$ for all $i$, so that $P(X_i=M_i)=0$ (resp., $>0$) for all $i$.
\end{Remark}

\begin{Example}[of \rem{example}]\rm 
Assume that $X_i\sim\exp(\lambda)$ for some $i$ and that $F_{X_j}(t)<1$ for all $t$ and $j\not=i$. Then necessarily, for $j\not=i$,
\begin{equation}
P(X_j>M_j)\ge P(X_j>X_i)=E\bar F_{X_j}(X_i)>0\,.
\end{equation} 
With $a^+=\max(a,0)$, $P(X_i>M_i)=Ee^{-\lambda M_i^+}>0$ and for all $j$, 
\begin{equation}
P(X_j<M_j)=E\prod_{k\not=j}\bar F_{X_k}(X_j)>0\,.
\end{equation} 
Thus, when in addition time and cause are independent, it necessarily follows that $X_1,\ldots,X_n$ must all be exponentially distributed. 

For the case $n=2$, if $X$ has an exponential distribution and $F_Y(0)<1$, then independence of time and cause is equivalent to $Y$ being exponentially distributed and there is no need for any further assumptions. To see this, observe that $P(X>Y)=Ee^{-\lambda Y^+}>0$ and by right continuity there exists some $t>0$ such that $F_Y(t)<1$. Therefore, \begin{equation}
P(X<Y)\ge P(X<t,t<Y)=(1-e^{-\lambda t})\bar F_Y(t)>0\,.
\end{equation}
\end{Example}

\begin{Remark}\rm
As in \rem{expgeo}, if $X_1,\ldots,X_n$ are independent and, for some strictly increasing $g$ and all $i$, $X_i\sim \gexp(\lambda_i)$  ($\lambda_i>0$) then, with $\lambda=\sum_{i=1}^n\lambda_i$, $M\sim \gexp(\lambda)$ and, for all $i$, $P(X_i=M)=P(X_i<M_i)=\lambda_i/\lambda$. Similarly, if we replace $\gexp(\lambda_i)$ by $\ggeo(\alpha_i)$ ($\alpha_i\in(0,1)$),
then, with $\alpha=1-\prod_{i=1}^n(1-\alpha_i)$ we have that $M\sim \ggeo(\alpha)$ and, for all $i$, 
\begin{align}
P(X_i<M_i)&=\frac{\alpha_i\prod_{j\not=i}(1-\alpha_j)}{\alpha}=\frac{\alpha_i(1-\alpha)}{(1-\alpha_i)\alpha}\nonumber\\
P(X_i=M_i)&=\frac{\alpha_i\left(1-\prod_{j\not=i}(1-\alpha_j)\right)}{\alpha}=\frac{\alpha_i(\alpha-\alpha_i)}{(1-\alpha_i)\alpha}\,,
\end{align} 
so that $P(X_i=M)=P(X_i\le M_i)=\alpha_i/\alpha$.

If we say that {\em $X_i$ is a cause} when $X_i=M$, then we see that the probability that $X_i$ is a cause is either $\lambda_i/\lambda$ or $\alpha_i/\alpha$, depending on the case and is independent of $g$.
\end{Remark}

\begin{Remark}\label{rem:1}\rm
When the pair $(X,Y)$ is exchangeable with $P(X=Y)=0$, $X\wedge Y$ is always independent of $\{Y>X\}$ (with $p=1/2$). In fact, in this case it is easy to check that $f(X,Y)$ is independent of $\{Y>X\}$ for any symmetric Borel $f$. 

In contrast to the independent case, the marginal distribution (the same for $X$ and $Y$) here need not be continuous even though $P(X=Y)=0$. For a trivial example, take 
\begin{equation}P(X=1,Y=0)=P(X=0,Y=1)=1/2\,,\end{equation} where the marginal distribution is $\text{Bernoulli}(1/2)$. 

In the i.i.d. case, the condition $P(X=Y)=0$ is equivalent to the continuity of their distribution, a case which is already contained in \thm{main} and \cite{s65}. We also note that other cases where $X,Y$ are dependent and have a joint absolutely continuous distribution were discussed in \cite{kp91}.
\end{Remark}

\begin{Remark}\rm 
If $F\sim \ggeo(\alpha)$, where $\alpha\in(0,1)$ and $g$ is strictly increasing, then, with $N(t)=\sup\{n|\ g(n)\le t\}$, we have that $F(t)=1-e^{-\lambda N(t)}$, where $\lambda=-\log(1-\alpha)$ (so $\alpha=1-e^{-\lambda}$) and thus $H_F(t)=\lambda N(t)$. This implies that \eq{CF} is valid for this case as well only that $V$ on the right is geometrically rather than exponentially distributed. Also note that if $X\sim \ggeo(1-e^{-\lambda})$ and $Y\sim \ggeo(1-e^{-\mu})$ then, with $t_F=\sup\{g(n)|\,n\ge 1\}$, we have that
$H_Y(t)/H_X(t)=\mu/\lambda$ for $t<t_F$, so that, like in the $g$-exponential case, the cumulative hazard functions are proportional and
$\bar F_Y(\cdot)=\bar F_X(\cdot)^{\mu/\lambda}$. However, unlike for the continuous case, where proportional cumulative hazard functions imply that $X,Y$ are necessarily $g$-exponentially distributed with the same strictly increasing $g$, it does not necessarily follow that $X$ and $Y$ are $g$-geometrically distributed when the distributions are not continuous or even completely atomic (that is, discrete).

Finally, recall that if $V\sim \exp(\lambda)$ then $\lceil V\rceil \sim \geo(1-e^{-\lambda})$ and thus the class of all $g$-geometric distributions is the same as the class of all $h$-exponential distributions for $h(t)=g(\lceil t\rceil)$. Of course, $h$ is not strictly increasing on $\mathbb{R}$.
\end{Remark}

\begin{Remark}\rm
Obviously, for independent $X,Y$ (similarly for the multivariate case), $\max(X,Y)$ is independent of $\{X>Y\}$, $\{X=Y\}$ and $\{X<Y\}$ iff independence of time and cause holds for $-X,-Y$. When $P(X>Y)>0$ and $P(X<Y)>0$, this happens iff $X,Y$ are $g$-exponentially or $g$-geometrically distributed, but now for for some strictly {\em decreasing} $g$. Is it possible for independence of time and cause to hold for both $X,Y$ and $-X,-Y$? Other than the trivial cases where either $P(X>Y)=1$ or $P(X<Y)=1$ or $P(X=Y)=1$, this is impossible.
\end{Remark}

\begin{Remark}\rm
Recalling $F_\lambda$ from \eq{CF}, some obvious families of $g$-exponential distributions are: 
\begin{equation}
 F_\lambda(x)=\begin{cases}
1-(1-x)^\lambda 1_{(0,1)}(x) &\text{Beta}(1,\lambda)\\ 
(1-e^{-\lambda x})1_{(0,\infty)}(x) &\exp(\lambda)\\
(1-e^{-\lambda x^\beta})1_{(0,\infty)}(x)&\text{Weibull}(\lambda,\beta)\\
1-\left(\frac{a}{x}\right)^\lambda 1_{(a,\infty)}(x)&\text{Pareto}(a,\lambda)\,.
\end{cases}
\end{equation}
where $a,\beta>0$. In these examples $g(x)=1-e^{-x},\,x,\,x^{1/\beta},\,ae^x$, respectively.
\end{Remark}

\begin{Remark}\rm
A tiny statistical aspect in passing. Recall \eq{finverse} and \eq{hinverse}. For a continuous $F$ we have that $F(F^{-1}(u))=u$ for all $u\in(0,1)$ (see~(iv)). Consequently, $H_F(H_F^{-1}(y))=y$ for all $y>0$. We saw that in this case $H_F^{-1}(V)\sim F$ where $V\sim\exp(1)$ and thus $V\sim H_F(X)$ where $X\sim F$. Therefore, if $X,Y$ are independent and $X\sim F$, then necessarily $P(X=Y)=0$ and then, provided that $P(X>Y)>0$ and $P(X<Y)>0$, time and cause are independent iff $Y\sim H^{-1}_F(W)$ for some exponentially distributed $W$ (recall that, by continuity of $F$, $H_F^{-1}$ is strictly increasing). Hence, $H_F(Y)\sim H_F(H^{-1}_F(W))=W$ is exponentially distributed. This means that for $X$ with a known continuous distribution $F$ and an independent random variable $Y$ such that $P(X>Y)>0$ and $P(X<Y)>0$, testing the hypothesis that time and cause are independent is equivalent to testing that $H_F(Y)$ is exponentially distributed (not necessarily with rate~$1$). Note that $F$ here, hence $H_F$, is considered known, while the cdf of $Y$ is considered unknown. Also note that we are not assuming a-priori that $Y$ has a continuous distribution, nor that $F$ has an inverse (whenever $F$ is strictly increasing on $\{x|\,F(x)\in(0,1)\}$ or, equivalently, $F^{-1}$ is continuous on $(0,1)$).
\end{Remark}

\end{document}